\documentstyle[11pt,psfig]{article}
 \newtheorem{theorem}{Theorem}[section]
 \newtheorem{lemma}[theorem]{Lemma}
 \newtheorem{corollary}[theorem]{Corollary}

\newtheorem{definition}[theorem]{Definition}
 \newenvironment{proof}{{\it Proof:\/}}{$\Box$\vskip 0.08in}

\newtheorem{conjecture}[theorem]{Conjecture}
\newtheorem{problem}[theorem]{Problem}

\newtheorem{proposition}[theorem]{Proposition}

\newcommand{\pct}[1]{}%{{\Psfig{figure=#1}}}
\newcommand{\Psfig}[1]{{\mbox{$\ \ $}}}

\begin{document}

 \begin{center} {\large\bf
Fundamentals of Kauffman bracket skein modules}
\end{center}
 \begin{center}
                      J\'ozef H.~Przytycki
\end{center}
 \vspace{0.2in}
\begin{quotation} 
{\bf Abstract}\\
Skein modules are the main objects of an algebraic topology based on knots
(or position). In the same spirit as  Leibniz we would call our 
approach {\it algebra situs}.
When looking at the panorama of skein modules\footnote
{This is an extended version of a part of the talk ``Panorama of skein
modules", given at {\it Low Dimensional Topology Conference}; 
Madeira, Portugal, January, 1998.}, 
we see, past the  rolling hills
of homologies and homotopies, distant mountains - the Kauffman bracket skein
module,  and farther off in the distance skein modules based on other 
quantum invariants. We concentrate here on the basic properties of 
the Kauffman bracket skein module; properties fundamental in further 
development of the theory.
In particular we consider the relative Kauffman bracket skein module,
and we analyze skein modules of $I$- bundles over surfaces.
\end{quotation}

\ \\
{\bf History of skein modules from my personal perspective}\\

I would like to use this opportunity, of informal presentation,\footnote{
I would like here to thank Hanna Nencka for a 
titanic task of organizing the Madeira's conference.}
to give my personal history of algebraic topology based on knots (a more
formal account was given in \cite{Pr-7}).

In July 1986 I left Poland invited by Dale Rolfsen for a visiting position
at UBC. In January of 1987, Jim Hoste gave a talk at the first Cascade Mountains Conference
(in Vancouver) and described his work on multivariable generalization of
the Jones-Conway ([HOMFLY]\cite{P-T}) polynomial of links in $S^3$. 
He was convinced
that his construction works for 2 colors when the first color is represented
only by a trivial component. He had already succeeded in the case of 
2-component 2-bridge links. His method, following Nakanishi, was to
analyze link diagrams in an annulus (the trivial component being $z$ axis).
We immediately noticed (with Jim) that the analogous construction for
the Kauffman bracket polynomial has an easy solution \cite{H-P-1}.
In March of 1987 I got a manuscript on ``Invariants of colored links" 
\cite{H-K}.
I read it carefully and I was trying to generalize the work of Hoste and 
Kidwell to $n+1$ colors (n of which are used to color a trivial link of 
n components).
This led my attention to the possible torsion\footnote{
In November of 1984, I gave a talk at Warsaw seminar about 
the recently discovered
Jones polynomial, and when describing the Jones skein relation and the
Alexander skein relation I was asked by Pawe{\l} Traczyk whether we really
need any restrictions on coefficients. I realized then that even if
restrictions are needed we should not assume them from the beginning but
instead we should analyze their character.}
related to the Dirac trick (more precisely with 2-torsion in the mapping 
class group of the 2-sphere). I discussed it with Darryl McCullough,
whom I visited in the first part of April. On the last weekend of my stay
in Oklahoma, I was struck by the idea that I am not really analyzing
colored links, but a (monochromatic) knot theory in the solid torus or
the connected sum of solid tori. Knots are formally added and taken
modulo the skein relation coming from the Jones-Conway polynomial, in 
the case of \cite{H-K}, or from 
the Kauffman bracket polynomial in the case of \cite{H-P-1}.
For me, this warm April day in Oklahoma was the birth 
of {\it skein modules}\footnote{I was hesitant about what to call
 these new objects. I thought that
Conway's ``linear skein" was to parochial for the ``big word" I envisioned
for the concept. Still I wanted to keep ``skein" acknowledging Conway's
vision \cite{Co-1}. 
The name ``skein group" would be natural (like homology group)
but misleading. Finally, I decided for the {\it skein module}.}.

Skein relations have their origin in an observation by Alexander 
(\cite{Al}, 1928) that his polynomials of 
three links $L_+, L_-$ and $L_0$ in $S^3$ are linearly related 
(here $L_+,L_-$ and $L_0$ denote three
links which are identical except in a small ball as shown
in Fig. 0.1). Conway rediscovered the Alexander observation and
normalized the Alexander polynomial so that it satisfies the skein relation
$$\Delta_{L_+}(z) - \Delta_{L_-}(z) = z\Delta_{L_0}(z)$$
 (\cite{Co-1}, 1969).
For my invention/discovery of skein modules, it was probably crucial that
I had read Conway's famous paper \cite{Co-1} and the following it work
by Giller \cite{Gi}, Kauffman \cite{Ka-1,Ka-2,Ka-3}, and 
Lickorish and Millett \cite{L-M-1} \footnote{Only later I learned that
in the late seventies Conway advocated the idea of
considering  the free $Z[z]$-module over oriented
links in an oriented 3-manifold and
dividing it by the submodule generated by his skein relation
\cite{Co-2} \cite{Co-3}; Conway called the resulting module ``linear skein".}.

\centerline{\psfig{figure=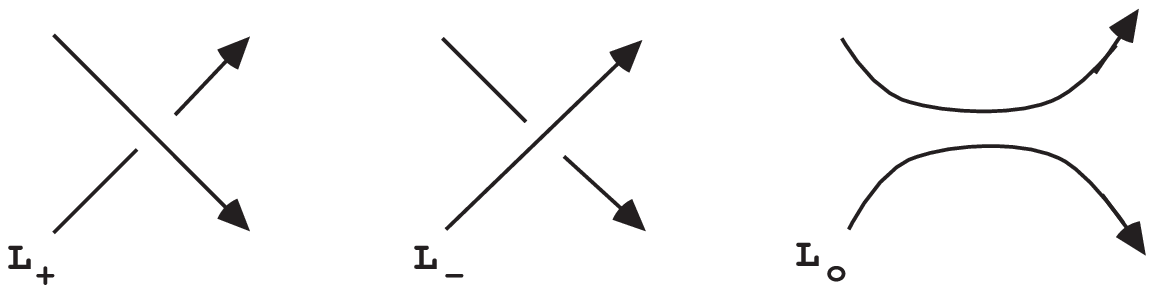}}
\begin{center}
Figure 0.1
\end{center}

Still in April of 1987 (visiting C.Gordon and J.Hempel in Texas), 
I interpreted several known facts in the language of skein modules. 
First was the observation that the existence of Jones type polynomials 
could be interpreted as saying that the related skein 
modules in $S^3$ were free and generated by the unknot.
Then the results of Hoste-Kidwell and Hoste-Przytycki interpret (compute)
 skein modules of the solid torus. 
The Kauffman method  allowed the computation of
the Kauffman bracket skein module of the product of a surface and 
the interval (this is very important fact which I will discuss in 
details later in the talk).
Relative  skein modules of the disc (Temperley-Lieb algebra,
Hecke algebra and Birman-Murakami-Wenzl algebra as they are known now) 
was computed/interpreted and shown to be free of 
$\frac{1}{n+1}{2n\choose n}$ (resp. $n!$ and $(2n-1)(2n-3)\cdot\ldots\cdot 1$) 
generators. I was told by Pawe{\l} Traczyk, in summer of 1986, about  
easy proofs of these facts.

I wrote the introductory paper on skein modules in May of 1987
\cite{Pr-1}. 
My initial definition of the Kauffman bracket skein module was rather
clumsy: it used unframed links up to regular isotopy. It worked well
for $M= F\times I$, in particular for a handlebody, but then I was
forced to consider a Heegaard decomposition of a manifold.

In May of 1988 (at the conference in Annapolis) I got a paper
by Vladimir Turaev, in which he introduced (independently) the concept
of skein modules \cite{Tu-1}. He pointed there importance of framing in
definitions of some skein modules.

\section{Skein modules of 3-manifolds}

Our goal is to build an algebraic topology based on knots
\footnote{One would like to say, in the spirit of Leibniz: {\it algebra
situs}. }. We call the main
object used in the theory {\it a skein module} and we associate it
to any 3-dimensional manifold.  
Skein modules are quotients of free modules over ambient
isotopy classes of links in a 3-manifold by properly chosen local
(skein) relations.
The choice of relations is a delicate task; we
should take into account several factors:
 
\begin{enumerate}
\item [(i)]  Is the module we obtain accessible (computable)?
\item [(ii)]
How precise are our modules in distinguishing 3-manifolds and links
in them?
\item [(iii)]
Does the module we obtain admit some additional structure (e.g.~filtration,
gradation, multiplication, Hopf algebra structure)?
\end{enumerate}
 
\noindent
From a practical point of view there is yet a fourth important factor
\begin{enumerate}
\item [(iv)]
The ``historical factor'' in the choice of (skein) relations:
the relations of links which were already studied (possibly for 
totally different reasons) will be compared with the new structures, 
just to see how they work in the new setup. 
For example, if we consider the Jones skein relation we
can be sure that even for $S^3$ we get a nontrivial result.
\end{enumerate}
 
The idea of the skein module should become more apparent after we consider
the main example of the talk, the Kauffman bracket skein module.
 
\section{The Kauffman bracket skein module}

The skein module based on the Kauffman bracket skein relation is, so far,
the most extensively studied object of the {\it algebraic topology based 
on knots}.
We describe in this section the basic properties of the Kauffman Bracket
Skein Module (KBSM) and list manifolds for which the structure of 
the module is known. In the third section, we give the detailed proof
of the structure of KBSM of a 3-manifold being an interval bundle over
a surface.
We extend the analysis to the case of the Relative Kauffman Bracket Skein
Module (RKBSM).
In the fourth section we discuss the torsion in KBSM. In particular,
we investigate in details the case of the RKBSM of a product of a surface
and the interval ($F\times I$).

\begin{definition}[\cite{Pr-1,H-P-3}]\label{2.1}\ \\
Let $M$ be an oriented 3-manifold, ${\cal L}_{fr}$ the set of unoriented
framed links in $M$ (including the empty knot, $\emptyset$), 
$R$ any commutative ring with identity and $A$ an
invertible element in $R$. Let $S_{2,\infty}$ be
the submodule of $R{\cal L}_{fr}$ generated by skein expressions
$L_+- AL_0 - A^{-1}L_{\infty}$, where the triple 
$L_+, L_0, L_{\infty}$ is presented in Fig.2.1, and $L \sqcup T_1 +(A^2 +
A^{-2})L$, where $T_1$ denotes the trivial framed knot. 
We define the Kauffman bracket skein module, ${\cal S}_{2,\infty}(M;R,A)$,  
as the quotient ${\cal S}_{2,\infty}(M;R,A)= R{\cal L}_{fr}/S_{2,\infty}$.
\end{definition}
\ \\
%\vspace*{0.8in}
\centerline{\psfig{figure=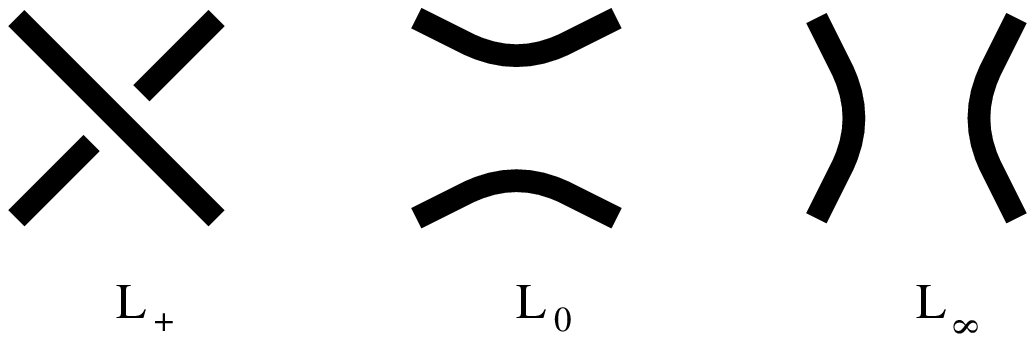}}
\begin{center}
Fig. 2.1.
\end{center}
Notice that $L^{(1)}=-A^3L$ in ${\cal S}_{2,\infty}(M;R,A)$,
where $L^{(1)}$ denotes a link obtained from $L$ by twisting the framing
of $L$ by a full twist in a positive direction. We call this the
framing relation. 
We use the simplified notation 
${\cal S}_{2,\infty}(M)$ for ${\cal S}_{2,\infty}(M;Z[A^{\pm 1}],A)$.

We list below several elementary properties of KBSM including
description of the 
KBSM of any compact 3-manifold using generators and relators. 
\begin{proposition} 
\begin{enumerate} 
\item [(1)] An orientation preserving embedding of 3-manifolds 
$i: M \to N$ yields the homomorphism of skein modules 
$i_*: {\cal S}_{2,\infty}(M;R,A) \to {\cal S}_{2,\infty}(N;R,A)$. 
The above correspondence leads to a functor from the category of 
3-manifolds and orientation preserving embeddings (up to ambient isotopy) 
to the category of $R$-modules (with a specified invertible element $A \in R$). 
\item [(2)] 
\begin{enumerate} 
\item [(i)] 
If $N$ is obtained from $M$ by adding a 3-handle to it (i.e. capping off 
a hole), and $i: M \to N$ is the associated embedding, 
then $i_*:{\cal S}_{2,\infty}(M;R,A) \to {\cal S}_{2,\infty}(N;R,A)$ 
is an isomorphism. 
\item [(ii)] If $N$ is obtained from $M$ by adding a 2-handle to it, 
and $i: M \to N$ is the associated embedding, 
then $i_*:{\cal S}_{2,\infty}(M;R,A) \to {\cal S}_{2,\infty}(N;R,A)$ 
is an epimorphism. 
\end{enumerate} 
\item [(3)] If $M_1 \sqcup M_2$ is the disjoint sum of 3-manifolds $M_1$ 
and $M_2$ then 
$${\cal S}_{2,\infty}(M_1\sqcup M_2;R,A)= {\cal S}_{2,\infty}(M_1;R,A) 
\otimes {\cal S}_{2,\infty}(M_2;R,A).$$ 
\item [(4)] 
(Universal Coefficient Property)\\ 
Let $r: R \to R'$ be a homomorphism of rings (commutative with 1). 
We can think of $R'$ as an $R$ module. Then the identity map on 
${\cal L}_{fr}$ induces the isomorphism of $R'$ (and $R$) modules: 
$$ \bar r: {\cal S}_{2,\infty}(M;R,A)\otimes_R R' \to 
{\cal S}_{2,\infty}(M;R',r(A))  .$$ 
\item [(5)] 
Let $(M,\partial M)$ be a 3-manifold with the boundary $\partial M$, 
and let $\gamma $ be a simple closed curve on the boundary. Let $N=M_{\gamma}$ 
be the 3-manifold obtained from $M$ by adding 
a 2-handle along $\gamma $. Furthermore let ${\cal L}_{fr}^{gen}$ be a set 
of framed links in $M$ generating ${\cal S}_{2,\infty}(M;R,A)$.\\ 
Then ${\cal S}_{2,\infty}(N;R,A) = {\cal S}_{2,\infty}(M;R,A)/J$, where 
$J$ is the submodule of ${\cal S}_{2,\infty}(M;R,A)$ generated by 
expressions $L- sl_{\gamma}(L)$, where $L\in {\cal L}_{fr}^{gen}$ and 
$sl_{\gamma}(L)$ is obtained from $L$ by sliding it along ${\gamma}$ 
(i.e. handle sliding). 
\item [(6)] 
Let $M$ be an oriented compact manifold and consider its Heegaard decomposition 
(that is $M$ is obtained from the handlebody $H_n$ by adding 2 and 3-handles 
to it), then $M$ has a presentation as follows: generators of 
${\cal S}_{2,\infty}(M;R,A)$ are generators of ${\cal S}_{2,\infty}(H_n;R,A)$ 
and relators are yielded by 2-handle slidings. 
\end{enumerate} 
\end{proposition} 
\begin{proof} 
\begin{enumerate} 
\item [(1)] $i_*$ is well defined because if framed links $L_1$ and $L_2$
are ambient isotopic in $M$ then $i(L_1)$ and $i(L_2)$ are ambient 
isotopic in $N$. Furthermore any skein triple $L_+,L_0,L_{\infty}$ in $M$, 
is sent by $i$ to a skein triple in $N$. Finally $i(T_1)$ is a trivial
framed knot in $N$.\ Notice that if $i_*:M \to N$ is an orientation reversing 
embedding then $i_*$ is a $Z$-homomorphism and $i(Aw)=A^{-1}i(w)$. 
\item [(2)] 
\begin{enumerate} 
\item [(i)] It holds because the cocore of a 3-handle is 
$0$-dimensional.\footnote{A manifold $N$ is obtained from an $n$-dimensional 
manifold $M$
by attaching a $p$-handle, $D^p \times D^{n-p}$, to $M$, if $N= M \cup_f
D^p \times D^{n-p}$ where $f: \partial D^p \times D^{n-p}$ is an embedding.
$D^p \times \{0\}$ is a core of the handle and 
$\{0\} \times D^{n-p}$ is a cocore of the handle \cite{R-S}.} 
\item [(ii)] It holds because the cocore of a 2-handle is $1$-dimensional. 
\end{enumerate} 
\item [(3)] This is a consequence of the well known property of short exact 
sequences, \cite{Bl}:\\ 
If $ 0 \to A' \to A \to A'' \to 0$ and $ 0 \to B' \to B \to B'' \to 0$ are 
short exact sequences of $R$-modules then 
$ 0 \to A'\otimes B + A\otimes B' \to A\otimes B \to A''\otimes B'' \to 0$ 
is a short exact sequence. 
\item [(4)] 
The exact sequence of $R$ modules 
$$ S_{2,\infty}(R,A) \to R{\cal L}_{fr} \to {\cal S}_{2,\infty}(M;R,A) \to 0$$ 
leads to the exact sequence of $R'$ modules (\cite {C-E}, Proposition 4.5): 
$$ S_{2,\infty}(R,A)\otimes_R R' \to R{\cal L}_{fr}\otimes_R R' 
\to {\cal S}_{2,\infty}(M;R,A) \otimes_R R' \to 0.$$ 
Now, applying the ``five lemma" to the commutative diagram with exact rows 
(see for example \cite{C-E} Proposition 1.1) 
\begin{displaymath} 
\begin{array}{ccccccc} 
 S_{2,\infty}(R,A) \otimes_R R' & \to & R{\cal L}_{fr}\otimes_R R' 
& \to & {\cal S}_{2,\infty}(M;R,A) \otimes_R R' & \to & 0 \\ 
\downarrow epi & & \downarrow iso & & \downarrow \bar r & & \\ 
S_{2,\infty}(R',\bar r (A)) & \to & R'{\cal L}_{fr} & \to & 
{\cal S}_{2,\infty}(M;R',\bar r(A))  & \to & 0 
\end{array} 
\end{displaymath} 
we conclude that $\bar r$ is an isomorphism of $R'$ (and $R$) modules. 
\item[(5)]
It follows from (2)(ii) because any skein relation can be performed in $M$,
and the only difference between KBSM of $M$ and $N$ lies in the fact that 
some nonequivalent links in $M$ can be equivalent in $N$; 
the difference lies exactly in the possibility 
of sliding a link in $M$ along the added 2-handle (that is $L$ is moving 
from one side of the cocore of the 2-handle to another).
\item[(6)] It follows from (5) and (2)(i).
\end{enumerate} 
\end{proof} 

In the next theorem we list manifolds for which the exact structure of
the Kauffman bracket skein module has been computed.
\begin{theorem}[\cite{Ka-4,Pr-1,H-P-4,H-P-5,H-P-6,Bu-2}]\label {2.3}\ \\
\begin{enumerate}
\item
[(a)]  ${\cal S}_{2,\infty}(S^3)=Z[A^{\pm 1}]$, 
more precisely: $\emptyset$ is the generator of the module
and $L=<L> T_1=(-A^2-A^{-2})<L> \emptyset$ where $<L>$ is the Kauffman 
bracket polynomial of a framed link $L$.
\item
[(b)]  ${\cal S}_{2,\infty}(F\times [0,1])$ is a free module generated
by links (simple closed curves) on  $F$ with no trivial component 
(but including the empty knot). 
Here $F$ denotes  an oriented surface (see also Theorems 3.1 and 3.9).\\ 
This applies in particular to a handlebody, because 
$H_n = P_n \times I$, where $H_n$ is a handlebody of genus $n$ and $P_n$ is 
a disc with $n$ holes. 
\item
[(c)]  ${\cal S}_{2,\infty}(L(p,q))$ is a free $Z[A^{\pm 1}]$ module and 
it has $[p/2]+1$ generators, where $[x]$ denotes the integer part of $x$. 
\item
[(d)] ${\cal S}_{2,\infty}(S^1\times S^2)=  Z[A^{\pm 1}]\oplus 
\bigoplus_{i=1}^{\infty} Z[A^{\pm 1}]/(1-A^{2i+4})$
\item
[(e)] The skein module of the complement of the torus knot
of type $(k,2)$ is a free $Z[A^{\pm 1}]$-module 
generated by links $u_{i,j}$ ($i\geq 0$, $\frac{k-1}{2}\geq j\geq0$ ), 
where $u_{i,j}$
is composed of $i$ meridians and $j$ curves $\gamma$, as illustrated
in Fig. 2.2. 

\item 
[(f)] Let $W$ be the classical Whitehead manifold, then ${\cal S}_{2,\infty}(W)$
is infinitely generated torsion free but not free.
\end{enumerate} 
\end{theorem}

\ \\
\centerline{\psfig{figure=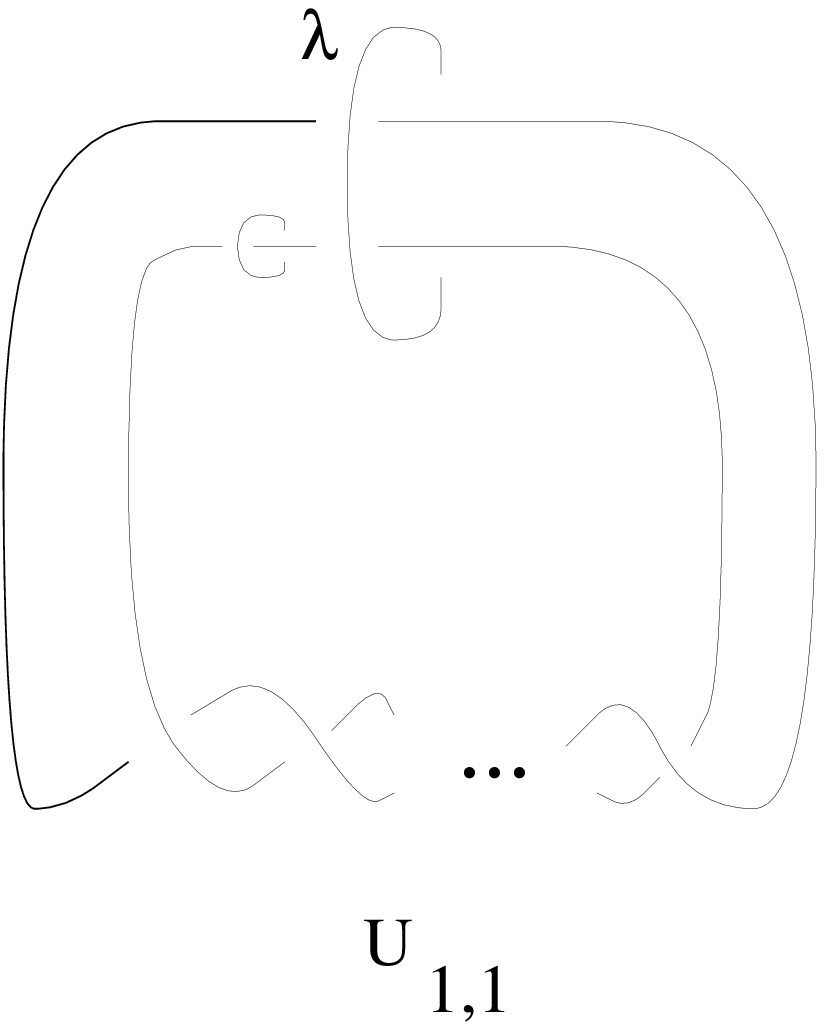,height=1.5in}}
%\vspace*{1.5in}
\begin{center}
Fig. 2.2.
\end{center}

In \cite{H-P-5}, (f) is proved for a large class of genus one Whitehead
type manifolds. For the classical Whitehead manifold it seems
feasible to find the exact structure of ${\cal S}_{2,\infty}(W)$, 
and we plan to address that in a future paper.

We prove (b), with its generalizations, in the next section.

%\newpage

\section{KBSM and relative KBSM of $F\times I$ and $F\hat{\times} I$ }\ \\

The understanding of the Kauffman bracket skein module of the product
of a surface and the interval is the first step to understanding KBSM
of a general 3-manifold. Furthermore the case of $F\times I$ is 
relatively easy to understand because we can project links onto the surface
and work with diagrams of links. This can be generalized to twisted 
$I$-bundles over $F$ and one can have reasonable hope that the method
can work for other 3-manifolds by projecting links to spines of 3-manifolds.
The relative case is described in Theorem 3.9.

\begin{theorem}\label{3.1}
Let $M$ be an oriented 3-manifold which is either equal to $F\times I$,
where $F$ is an oriented surface, or it equal to a twisted $I$ bundle over
$F$ ($F\hat{\times} I$), where $F$ is an unoriented surface. Then the KBSM, 
${\cal S}_{2,\infty}(M;R,A)$, is a free $R$-module with 
a basis $B(F)$ consisting of links in $F$ without contractible 
components (but including the empty knot). 
\end{theorem}

\begin{proof} We will give here the proof of Theorem 3.1 which is based on the
original proof of Kauffman on the existence of his bracket polynomial.
Let $M$ be an oriented 3-manifold which is an $I$-bundle over a 
surface $F$\footnote{Because $M$ is oriented therefore 
for $\gamma$ in $F$ changing 
orientation of $F$, the restriction of the $I$-bundle to $\gamma$ is a
nontrivial bundle (M\"obius band). For $\gamma$ preserving orientation of $F$,
the bundle is trivial (an annulus).}.
Let $B(F)$ consist of all links in $F$ which have no trivial components
(including $\emptyset$). Furthermore each link is equipped with
an arbitrary, but specific framing (to be concrete we can assume that if
a knot in $F$ preserves the orientation of $F$ then we choose as its framing
the regular neighborhood of $K$ in $F$ (``blackboard" framing), if $K$ is
changing the orientation on $F$ than its regular neighborhood is a M\"obius
band so to get a framing we perform a positive half twist on it). 
Now one can quickly see that $B(F)$ is a generating set of 
${\cal S}_{2,\infty}(M;R,A)$. 
Namely every link
in $M$ has a regular projection on $F$ and any link can be reduced by  skein 
relations so that a projection has no crossings. 
Then another relation allows us to eliminate trivial components and finally 
the framing relations allow us to adjust framing. 
We will prove that $B(F)$ is a basis for 
${\cal S}_{2,\infty}(M;R,A)$. First we need however to consider the space of 
link diagrams (for a nonorientable surface $F$ the proof is still 
straightforward but requires great care).
\begin{definition}\label{3.2}
\begin{enumerate}
\item[(a)]
A link diagram (or marked diagram) on $F$ is a 4-valent graph in $F$ 
(allowing loops without vertices) such that one corner of a neighborhood of
each vertex 
is marked. $F$ does not need to be oriented for this definition. 
\item[(b)] 
Let $\cal D$ be a set of link diagrams on $F$ (up to isotopy of $F$),
and $R\cal D$ the free module over $\cal D$. The skein space of
diagrams, ${\cal S}\cal D$ is defined as a quotient:
$${\cal S}{\cal D}(F;R,A)= R{\cal D}/
({\psfig{figure=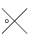,height=0.4cm}} -
A{\psfig{figure=L0nmaly.eps}} - A^{-1} {\psfig{figure=Linftynmaly.eps}},\ \
 D\sqcup T_1 +(A^2 + A^{-2})D).$$
\end{enumerate}
\end{definition}

\begin{lemma}\label{3.3}\ \\
Let $B^d(F)$ denote the set of link diagrams in $F$ without vertices 
and without trivial components (but allowing $\emptyset$). We can
identify $B^d(F)$ with the set $B(F)$ with framing ignored.
$B^d(F)$ is a subset of the set of link diagrams, so we have a homomorphism 
$\phi: RB^d(F) \to {\cal S}{\cal D}(F;R,A)$ defined by
associating to a link diagram in $F$, $\gamma \in B^d(F)$, its class in 
${\cal S}{\cal D}(F;R,A)$.\\
Then $\phi$ is an isomorphism.
\end{lemma}

\begin{proof} For any $D \in \cal D$ we can use the first relation to 
eliminate all crossings, and the second to eliminate trivial components of
$D$. Thus $\phi$ is an epimorphism.\\
To show that it is a monomorphism we will construct the inverse 
map, $\psi$.\\
First we define a map $\hat\psi: R{\cal D} \to RB^d(F)$.
Let $D\in \cal D$. We define $\hat\psi(D)$ as follows:\\
Choose any ordering $p_1,...,p_n$ of crossings of $D$, and use the
formula $D^{p_i}_{{\psfig{figure=Lmarkednmaly.eps,height=0.4cm}}}=
AD^{p_i}_{{\psfig{figure=L0nmaly.eps}}} + 
A^{-1} D^{p_i}_{{\psfig{figure=Linftynmaly.eps}}}$, for each crossing,
until all crossings are eliminated.
The upper index denotes the crossing
at which we perform a smoothing (crossing elimination). 
The result does not depend on the order of the crossings since we can make
any transposition of adjacent (with respect to ordering), pairs 
and get the same result:\\
\ \\
$(D^p_{{\psfig{figure=Lmarkednmaly.eps,height=0.4cm}}})^
{q}_{{\psfig{figure=Lmarkednmaly.eps,height=0.4cm}}} =
A(D^p_{{\psfig{figure=L0nmaly.eps}}})^
q_{{\psfig{figure=Lmarkednmaly.eps,height=0.4cm}}} +
A^{-1}(D^p_{{\psfig{figure=Linftynmaly.eps}}})^
q_{{\psfig{figure=Lmarkednmaly.eps,height=0.4cm}}}
=\\
\ \\
A^2(D^p_{{\psfig{figure=L0nmaly.eps}}})^q_{{\psfig{figure=L0nmaly.eps}}} +
(D^p_{{\psfig{figure=L0nmaly.eps}}})^q_{{\psfig{figure=Linftynmaly.eps}}} +
(D^p_{{\psfig{figure=Linftynmaly.eps}}})^q_{{\psfig{figure=L0nmaly.eps}}} +
A^{-2}(D^p_{{\psfig{figure=Linftynmaly.eps}}})^q_
{{\psfig{figure=Linftynmaly.eps}}}$\\ 
\ \\
and \\
$(D^q_{{\psfig{figure=Lmarkednmaly.eps,height=0.4cm}}})^
{p}_{{\psfig{figure=Lmarkednmaly.eps,height=0.4cm}}} =   
A(D^q_{{\psfig{figure=L0nmaly.eps}}})^
p_{{\psfig{figure=Lmarkednmaly.eps,height=0.4cm}}} + 
A^{-1}(D^q_{{\psfig{figure=Linftynmaly.eps}}})^
p_{{\psfig{figure=Lmarkednmaly.eps,height=0.4cm}}}  
=\\
\ \\                                                                      
A^2(D^q_{{\psfig{figure=L0nmaly.eps}}})^p_{{\psfig{figure=L0nmaly.eps}}} + 
(D^q_{{\psfig{figure=L0nmaly.eps}}})^p_{{\psfig{figure=Linftynmaly.eps}}} + 
(D^q_{{\psfig{figure=Linftynmaly.eps}}})^p_{{\psfig{figure=L0nmaly.eps}}} + 
A^{-2}(D^q_{{\psfig{figure=Linftynmaly.eps}}})^p_ 
{{\psfig{figure=Linftynmaly.eps}}}$\\ 
\ \\
After smoothing all crossings we eliminate trivial components by the relation
 $D\sqcup T_1 = (-A^2 -A^{-2})D$ (there is no ambiguity in the reduction).
Thus $D$ is uniquely expressed as a linear combination of elements of
$B^d(F)$, and we define $\hat\psi(D)$ as this linear combination (which 
lies in $RB^d(F)$). Therefore $\hat\psi$ is well defined.
Now $\hat\psi$ descends to $\psi: {\cal S}{\cal D}(F;R,A) \to RB^d(F)$ because
$\hat\psi ( {\psfig{figure=Lmarkednmaly.eps,height=0.4cm}} - 
A{\psfig{figure=L0nmaly.eps}} - A^{-1} {\psfig{figure=Linftynmaly.eps}})=0$
and $\hat\psi ( D\sqcup T_1 +(A^2 + A^{-2})D)=0$. 
Now, obviously, $\psi\phi = Id$, thus $\phi$ is a monomorphism.
\end{proof}
  Our goal is to prove that $B(F)$ is a basis of the Kauffman bracket
skein module  ${\cal S}_{2,\infty}(M;R,A)$, where $M$ is an oriented
3-manifold being an $I$ bundle over a surface $F$. Because we would like to
consider the case of orientable and unorientable surface simultaneously,
it is convenient to consider half-integer framings of links, that
is, to allow embedded M\"obius bands. This suggests the following
definition.
\begin{definition}\label{3.4}
Let $M$ be any oriented 3-manifold, $\bar{\cal L}_{fr}$ the set of embeddings
of annuli and M\"obius bands in $M$ (up to an ambient isotopy of $M$) 
and $\bar R$ 
a commutative ring with identity with a chosen invertible element $\bar A$ 
(we define $A=-\bar A^2$ and we will often write $\sqrt{-A}$ for $\bar A$. 
Let $\bar R \bar{\cal L}_{fr}$ denote a free $\bar R$ module over 
$\bar{\cal L}_{fr}$ and let $\bar S_{2,\infty}$ denote the submodule of
$\bar R\bar{\cal L}_{fr}$ generated by expressions 
$L_+ -AL_- - A^{-1}L_{\infty}$,
and $L^{1/2} - (\sqrt{-A})^3L$, where $L^{1/2}$ denotes $L$ with its 
framing twisted by a half twist in a positive direction. 
As before, for convenience, we allow the empty knot, $\emptyset$, and
add the relation $T_1=(-A^2-A^{-2})\emptyset$. \\
Then we define $\bar{\cal S}_{2,\infty}(M,\bar R,\bar A)= 
\bar R \bar{\cal L}_{fr}/\bar S_{2,\infty}.$
\end{definition}

Consider the $\bar R$-homomorphism
$g: {\cal S}{\cal D}(F;\bar R, A) \to \bar{\cal S}_{2,\infty}(M,\bar R,\bar A)$
defined on the basic elements $\gamma \in B^d(F)$ by $g(\gamma)=
\gamma^{fr}$ where $\gamma^{fr}$ is a framed link obtained from
$\gamma$ by giving it the blackboard framing (it may be an annulus or a 
M\"obius band). Using our skein relations,
in a similar manner as before, we see that $g$ is an epimorphism.
If $D$ is any marked diagram we can describe the framed link $g(D)$ as
follows: we resolve every crossing of $D$ according to the rule given
in Fig. 3.1 and giving the link $g(D)$ the blackboard framing
(the orientation of $M$ in a neighborhood of the crossing should
agree with that of $R^3$ from Fig. 3.1).\ \\
\ \\
\centerline{\psfig{figure=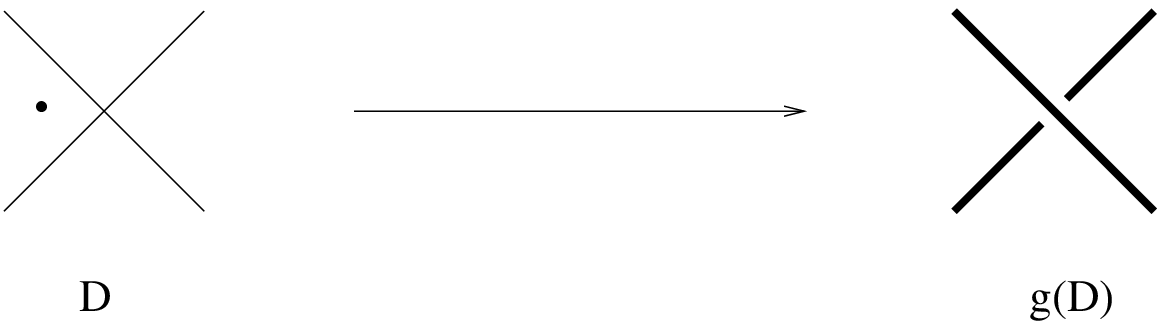}}
%\vspace*{1.2in}
\begin{center}
Fig. 3.1.
\end{center}

Consider now the following lemma
concerning Reidemeister moves on diagrams.
\begin{lemma}\label{3.5}
Consider the moves $\bar R_1, \bar R_2, \bar R_3$ on marked diagrams
(described below).
In ${\cal S}{\cal D}(F;\bar R, A)$ they satisfy:
\begin{enumerate}
\item[($\bar R_1$)]\ {\psfig{figure=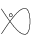}} 
$= - A^3$ {\psfig{figure=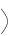}}
and {\psfig{figure=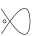}} 
$=- A^{-3}$ {\psfig{figure=Lline.eps}}, where
$\bar R_1$({\psfig{figure=Lline.eps}})$=$ {\psfig{figure=Ltwistgora.eps}} or
{\psfig{figure=Ltwistdol.eps}}. 
\item[($\bar R_2$)]\ $\bar R_2(D)=D$, where 
$\bar R_2$({\psfig{figure=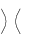}})
$=$
{\psfig{figure=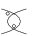}}.
\item[($\bar R_3$)]\  $\bar R_3(D)=D$, where 
$\bar R_3$({\psfig{figure=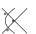}})
$=$
{\psfig{figure=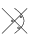}}.
\end{enumerate} 
\end{lemma}
\begin{proof}
\begin{enumerate} 
\item[($\bar R_1$)]\ \\ 
{\psfig{figure=Ltwistgora.eps}} $= A$\ {\psfig{figure=Lline.eps}}\ 
$\sqcup O +A^{-1}${\psfig{figure=Lline.eps}} $= (A(-A^2-A^{-2}) + A^{-1})$
{\psfig{figure=Lline.eps}} $= - A^3$ {\psfig{figure=Lline.eps}}.\\
{\psfig{figure=Ltwistdol.eps}} $= A^{-1}$ {\psfig{figure=Lline.eps}}\ 
$\sqcup O + A${\psfig{figure=Lline.eps}} $= (A+A^{-1}(-A^2-A^{-2}))$ 
{\psfig{figure=Lline.eps}} $=- A^{-3}$ {\psfig{figure=Lline.eps}}.
\item[($\bar R_2$)]\ \\
{\psfig{figure=R2marked.eps,height=0.4cm}}
$=A$ {\psfig{figure=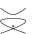,height=0.4cm}} $+ A^{-1}$
{\psfig{figure=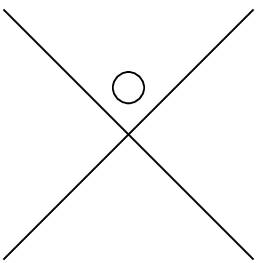,height=0.4cm}} 
$= A(A$ {\psfig{figure=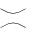,height=0.4cm}} $+ A^{-1}$
{\psfig{figure=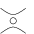,height=0.4cm}}) 
$+A^{-1}(A$ {\psfig{figure=L2lines.eps,height=0.4cm}}
$+A^{-1}$ {\psfig{figure=L2linesh.eps,height=0.4cm}}) 
$= (A^2 + AA^{-1}(-A^2 -A^{-2})+A^{-2})$ 
{\psfig{figure=L2linesh.eps,height=0.4cm}} $+$
{\psfig{figure=L2lines.eps,height=0.4cm}}
$=${\psfig{figure=L2lines.eps,height=0.4cm}}.
\item[($\bar R_3$)]\ \\
{\psfig{figure=R3markedb.eps,height=0.4cm}}$= A$
{\psfig{figure=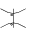,height=0.4cm}}
$+ A^{-1}$
{\psfig{figure=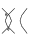,height=0.4cm}} 
$=A$
{\psfig{figure=R3markedbAsplit.eps,height=0.4cm}}
$+ A^{-1}$
{\psfig{figure=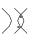,height=0.4cm}}
$=$ 
{\psfig{figure=R3markede.eps,height=0.4cm}}. 
We use here the invariance under $\bar R_2$
moves.
\end{enumerate}
\end{proof}
To use the lemma in the proof that $g$ is a monomorphism, we need a
variant of Reidemeister's theorem for marked diagrams:

\begin{proposition}\label{3.6}
Let $\hat g: {\cal D} \to \bar{\cal L}_{fr} $ be a map given by Fig. 3.1.
Then two marked diagrams, $D_1$ and $D_2$, represent 
the same framed link, $\hat g(D_1)=\hat g(D_2)$, \\ if and only if\\
one can go from $D_1$ to $D_2$ using Reidemeister moves $\bar R_i^{\pm 1}$ 
and an isotopy of $F$, and additionally, for corresponding link components
of $D_1$ and $D_2$, their Tait numbers are the same. One should notice here
that for a knot diagram the Tait number is independent on orientation
of the knot. Precisely for a knot diagram $D$ we define 
Tait($D$)$= \Sigma_p sgn(p)$, where 
sgn({\psfig{figure=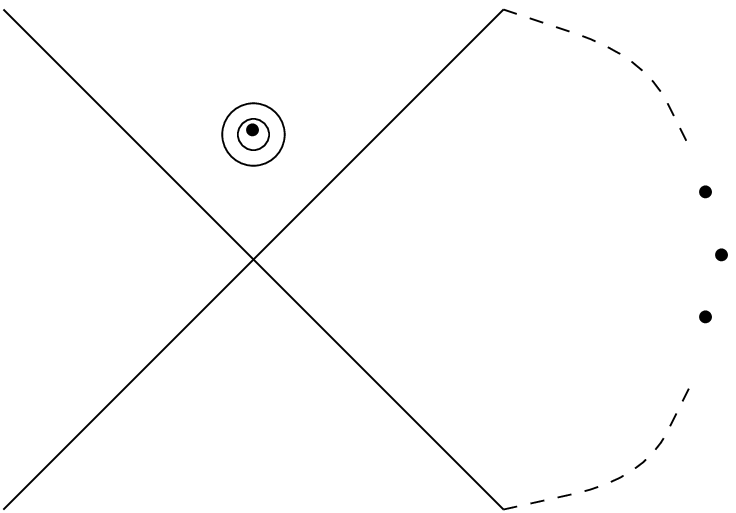,height=0.4cm}})$=1$ and 
sgn({\psfig{figure=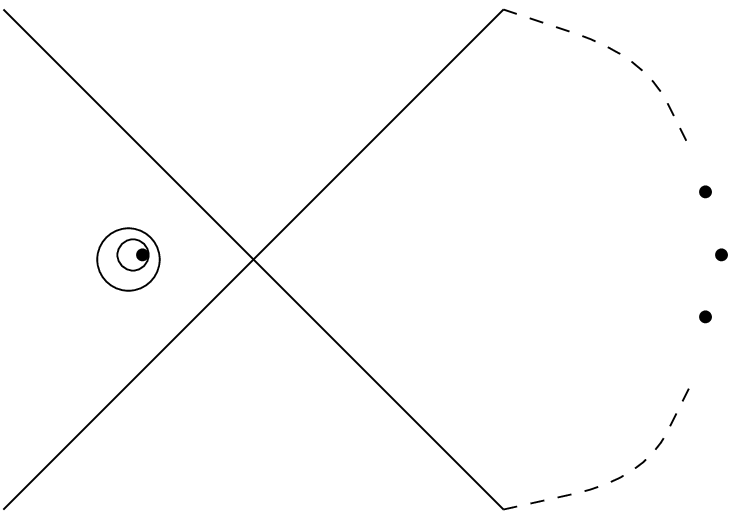,height=0.4cm}})$=-1$.

\end{proposition}
%\end{proof}
%\end{document}

\begin{proof} The proposition can be deduced from the classical Reidemeister
theorem and the result from the PL topology; Theorem 6.2 in
\cite{Hud}\footnote{
It follows from the theorem that if $C$ is a compact subset of a manifold 
$M$ and $F:M\times I \to M$ is the isotopy of $M$ then there is another
isotopy $\hat F:M\times I \to M$ such that\\
 $F_0 ={\hat F}_0$,
$F_1/C = {\hat F}_1/C$ and there exists a number $N$ such that
the set\ $\{x\in M \ |\ \hat F/\{x\}\times (k/N,(k+1)/N) \ is \ 
not \ constant\}$
sits in a ball embedded in $M$; \cite{Hud}, Corollary 6.3.}
\end{proof}
Our goal is to show that the epimorphism
$g: {\cal S}{\cal D}(F;\bar R, A) \to \bar{\cal S}_{2,\infty}(M,\bar R,\bar A)$
is an isomorphism.
We use Lemma 3.5 and Proposition 3.6 to construct the map inverse to
$g$.\
Let $\hat h: \bar R \bar{\cal L}_{fr} \to {\cal S}{\cal D}(F;\bar R,A)$ be a 
homomorphism defined as follows: choose a representative of a link 
$L\in \bar{\cal L}_{fr}$ which has a regular projection on $F$. Let $D_L$ be
a marked diagram on $F$ constructed as in Fig. 3.1, and let $t(L)$ be
the number (possibly half-integer) of positive twists which should 
be performed on the blackboard
framing of $D_L$ to get the framing of $L$. Then we define
$\hat h(L)=(-A^3)^{t(L)}D_L$. $\hat h(L)$ is well defined by Lemma 3.5 
and Proposition 3.6.\\
Furthermore $\hat h(L_+ -AL_- -A^{-1}L_{\infty})=0$, 
$\hat h(L \sqcup T_1 +(A^2 + A^{-2})L)=0$ and 
$\hat h(L^{1/2}-{\sqrt {-A^3}}L)=0$ so $\hat h$ descends to
$h:{\cal S}_{2,\infty}(M;R,A) \to {\cal S}{\cal D}(F;\bar R,A)$.
Of course $hg=Id$ so $g$ is a monomorphism, as required.
\end{proof}
We can now complete the proof of Theorem 3.1. Because 
$g: {\cal S}{\cal D}(F;\bar R, A) \to \bar{\cal S}_{2,\infty}(M,\bar R,\bar A)$
is an isomorphism, therefore $g(B^d(F))$ is a basis of 
$\bar{\cal S}_{2,\infty}(M,\bar R,\bar A)$. On the other hand $B(F)$, whose
elements may differ from elements of $g(B^d(F))$ only by framing, also forms
a basis of $\bar{\cal S}_{2,\infty}(M;\bar R,\bar A)$. Thus they are linearly
independent in ${\cal S}_{2,\infty}(M;R,A)$. Because $B(F)$ generates
${\cal S}_{2,\infty}(M;R,A)$ it is a basis of this module. The proof of
Theorem 3.1 is completed.

As an immediate corollary of Theorem 3.1 we obtain the structure of
 KBSM of the projective space $RP^3$. 
\begin{corollary}\label{3.7} 
${\cal S}_{2,\infty}(RP^3;R,A)=R\oplus R$. As a basis
of KBSM we can take $\emptyset$ and a generator of the fundamental group
of $RP^3$.
\end{corollary}
\begin{proof}
By Proposition 2.2(i) ${\cal S}_{2,\infty}(RP^3;R,A)=
{\cal S}_{2,\infty}(RP^3-int(D^3);R,A)$ and $RP^3-int(D^3)$ is equal to 
the twisted $I$-bundle over a projective plane ($RP^2\hat I$. By
Theorem 3.1, ${\cal S}_{2,\infty}(RP^2\hat I;R,A)$ is a free $R$-module
with basis $B(RP^2)$, which has two elements: the empty knot and the
noncontractible curve on $RP^2$. 
\end{proof}
One can generalize Theorem 3.1 to relative skein modules, as long
as a surface $F$ has a boundary. 
%boundary points of relative links are on the same boundary
%component of a manifold.

\begin{definition}[Relative Kauffman Bracket Skein Module]\ \\
Let $x_1,x_2,...,x_{2n}$ be a set of $2n$ (framed\footnote{
A framed point in $\partial M$ is an interval in $\partial M$. Thus
a relative framed link intersects $\partial M$ in framed points.}) 
points in $\partial M$,
where $M$ is an oriented 3-manifold. Let ${\cal L}_{fr}(n)$ be a family
of relative framed links in $(M,\partial M)$ such that 
$L \cap \partial M = \partial L =\{x_i\}$, 
considered up to an ambient isotopy fixing $\partial M$.
Let $R$ be a commutative ring with identity and $A$ an 
invertible element in $R$. Let $S_{2,\infty}(n)$ be 
the submodule of $R{\cal L}_{fr}(n)$ generated by the Kauffman bracket
skein relations. We define the Relative Kauffman Bracket Skein Module (RKBSM)
as the quotient:
$${\cal S}_{2,\infty}(M,\{x_i\}_1^{2n};R,A)= 
R{\cal L}_{fr}(n)/S_{2,\infty}(n)$$

\end{definition}

We list below a few useful properties of relative skein modules:
\begin{proposition}
\begin{enumerate}
\item[(a)] There is a functor from the category of
oriented 3-manifolds with $2n$ framed points on the boundary and 
orientation preserving embeddings (up to ambient isotopy fixed on the boundary)
to the category of $R$-modules (with a specified invertible element $A \in R$).
The functor sends an embedding $i: (M,\{x_i\}_{i=1}^{2n}) \to
(N,\{y_i\}_{i=1}^{2n})$ into $R$-modules morphism
${\cal S}_{2,\infty}(M,\{x_i\}_1^{2n};R,A) \to
{\cal S}_{2,\infty}(N,\{y_i\}_1^{2n};R,A)$.
\item[(b)] Adding a 3-handle to $M$ (outside $x_i$) is not changing the RKBSM,
and adding a 2-handle is adding only relations to RKBSM (handle slidings 
yield relations); compare Proposition 2.2(2). 
\item[(c)] The relative KBSM depends only on the distribution of boundary
points $\{x_i\}$ among boundary components of $M$, but not on the exact
position of $\{x_i\}$. In particular if $\partial M$ is connected, we can
write shortly ${\cal S}_{2,\infty}(M,n;R,A)$ instead of 
${\cal S}_{2,\infty}(M,\{x_i\}_1^{2n};R,A)$
\item[(d)] The relative KBSM satisfies Universal Coefficient Property,
compare Proposition 2.2(4).
\item[(e)] For a disjoint sum of 3-manifolds we have:
$${\cal S}_{2,\infty}(M_1\sqcup M_2,\{x_i,y_i\}_1^{2n};R,A)= 
{\cal S}_{2,\infty}(M_1,\{x_i\}_1^{2n};R,A) 
\otimes {\cal S}_{2,\infty}(M_2,\{y_i\}_1^{2n};R,A).$$
\end{enumerate}
\end{proposition}

\begin{theorem}\label{3.9} Let $M= F\bar{\times} I$ that is
 $M= F\times I$ or $M= F\hat{\times} I$, then
\begin{enumerate}
%\item [(a)] If all points $x_i$ lie on the same component of $\partial M$ then 
%${\cal S}_{2,\infty}(M,\{x_i\}_1^{2n};R,A)$ is a free $R$-module.
\item [(a)] Let $\partial F \neq \emptyset$ then
${\cal S}_{2,\infty}(M,\{x_i\}_1^{2n};R,A)$ is a free $R$-module. 
Consider all $x_i$ to
lie on $\partial F \times \{\frac{1}{2}\}$ then the basis of the
module ${\cal S}_{2,\infty}(M,\{x_i\}_1^{2n};R,A)$ is composed of
relative links on $F$ without trivial components.
\item [(b)] In the case of $F_{g,0}$ closed surface of genus $g$
 ($F\neq S^2$) the situation is more
delicate so we stop on the following observation:\\
${\cal S}_{2,\infty}(F_{g,0}\bar{\times} I;\{x_i\}_1^{2n};R,A)=
{\cal S}_{2,\infty}(F_{g,1}\bar{\times} I;\{x_i\}_1^{2n};R,A)/(I)$ where
$F_{g,1}=F_{g,0}-int(D^2)$ and assuming $x_i \in \partial D^2$, ideal $(I)$
is generated by moves in which arcs go above $D^2$.
%Even for $F=T^2$ and $n=1$ proof that the module is free 
%is rather involved but similar to Kauffman
% original proof.
\end {enumerate}
\end{theorem}
\begin{proof} The proof of (a) is the same as that of Theorem 3.1;
as before relative link diagrams representing the same link are
related by Reidemeister moves. In the case (b) it is no longer true as
we need also handle sliding. $F_{g,0}\bar{\times} I$ is obtained
from $F_{g,1}\bar{\times} I$ by adding the 2-handle along $\partial D^2$.
Now (b) follows from Proposition 3.8(b).  
\end{proof}
In the case $F$ is a closed surface the question whether 
${\cal S}_{2,\infty}(F \bar \times I ,\{x_i\}_1^{2n};R,A)$ is free is 
open in general. If not all $x_i$ lie on the same boundary component
of $F \times I $ then the skein module has a torsion in the case
of $F$ being a sphere or a torus. We prove it in the last section.
 We propose the following conjecture, which easily holds for $F=S^2$ and
otherwise
which we are able to confirm only for $F$ being a torus and $n=1$.
\begin{conjecture}\label{3.11}
Let $F$ be a closed surface and $x_i\in F\times \{0\}$ for any $i$, then
the skein module ${\cal S}_{2,\infty}(F \times I ,\{x_i\}_1^{2n};R,A)$ 
is free.
\end{conjecture}

Theorem 3.10 can be nicely illustrated by the following corollary.
\begin{corollary}\label{3.12}
\begin{enumerate} 
\item [(a)] ${\cal S}_{2,\infty}(D^2\times I,\{x_i\}_1^{2n};R,A)$ is
a free $R$ module of $\frac{1}{n+1} {{2n}\choose{n}}$ free generators.

\item [(b)] 
${\cal S}_{2,\infty}({\psfig{figure=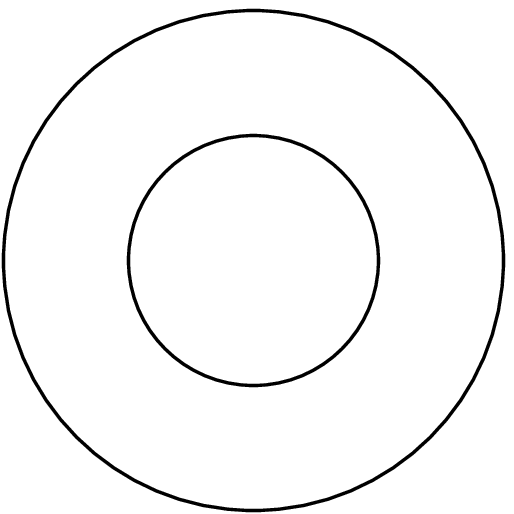,height=0.4cm}}
\ \times I,\{x_i\}_1^{2n};R,A)$ is
a free $R[\alpha]$ module with ${2n}\choose{n}$ free generators,
where $\alpha$ is represented by a longitude of the annulus.

\end{enumerate}
\end{corollary}
\begin{proof}
Corollary 3.12 (a) describes, well known module structure of the 
Temperley-Lieb algebra with basis consisting of the nth Catalan number 
of elements.\
(b) follows from the work of Jones and Tom Dieck [Jo,To]. 
We provide here a short, self-contained proof. 
In lieu of Theorem 3.9,
it suffices to count the crossless connections (by arcs) of $2n$ points 
in the disc and annulus. We offer an amazingly simple calculation 
for both cases simultaneously. 
Let $C_n$ be the number of connections in the disc and
$D_n$ be the number of connections in the annulus (all points $x_i$ are
on the "outside" circle of the annulus).
Connection arcs in the disc cut the disk into $n+1$ pieces.
To get a connection in
the annulus we have to put a ``table" in $D^2$ (remove a disk from $D^2$).
Thus $D_n = (n+1)C_n$. On the other hand, any arc of a connection in 
the annulus has a first point (with respect to some fixed orientation of 
the annulus) and any choice of $n$ points leads to a unique connection system,
for which given points are first.\footnote{The easy way to see this 
unique system of arcs which comes from a choice 
of $n$ points is to imagine the following childrens game: $n$ girls and $n$
 boys stand by the wall of a room. The result of the game is
that each girl will give a hand to a boy to her right (no crossings), but
not necessarily to the boy on her immediate right. 
In the first round each girl looks to her
right and gives her hand, if she has
a boy to her immediate right. Connected pairs are out of the game
which restarts with the remaining children.
At the end of the game all of the children are
paired up (exactly as needed in a connection system). No table was needed
in this game (algorithm). The table was necessary for the inverse function:
given connection arcs, decide on which end is a girl and on which end is 
a boy.}
Therefore $D_n= {2n \choose n}$ and thus $C_n= \frac{1}{n+1} {{2n}\choose{n}}$.
\end{proof}

For an oriented  surface $F$, let $x_1$, $x_2$,\ldots , $x_n$ 
$\in F\times \{1\}$ and $x_{n+1}$, $x_{n+2}$,\ldots , $x_{2n}$ 
$\in F\times \{0\}$, where $x_i$ and $x_{n+i}$ project to the same point
of $F$. Then ${\cal S}_{2,\infty}(F\times I,\{x_i\}_{i=1}^{2n};R,A)$ 
has an algebra structure
with the product of relative links $L_1 \cdot L_2$ defined by placing
$L_1$ in $F\times [\frac{1}{2},1]$ and $L_2$ in $F\times [0,\frac{1}{2}]$.
For $F=D^2$ the algebra is the Temperley-Lieb algebra, $TL_n$. For
$F= {\psfig{figure=annulus.eps,height=0.3cm}}$ we call the algebra,
annular Temperley-Lieb algebra and denote by $ATL_n$. The algebra related to
$ATL_n$ was first computed by Jones in the context of affine
Hecke algebras \cite{Jo} and independently by Tom Dieck (who denoted it by 
$TB_n$), \cite{To}. Corollary 3.12 supports a relatively short proof of
the structure of $ATL_n$. 

\begin{theorem}[{[Jo,To]}]\label{3.13}\ \\
\begin{enumerate}
\item [(a)]
 $ATL_n = {\cal S}_{2,\infty}({\psfig{figure=annulus.eps,height=0.3cm}}\times
I,\{x_i\}_{i=1}^{2n};R,A)=R[\alpha]<t,e_1,...e_{n-1}>/{\cal I}(n)$ 
where ${\cal I}(n)$ is the ideal generated
by the expressions:\\
$e_ie_j=e_je_i$ for $|j-i|\geq 2$, $e_ie_je_i=e_i$
for $|j-i|=1$, $te_i=e_it$ \\ for $i>1$,  $e_1te_1=\alpha e_1$,
 $e_i^2=(-A^2-A^{-2})e_i$,
$t^2=-A^{-2}\alpha t-A^{-4}$.
\item [(b)] $ATL_n$, as a module, is freely generated by 
${2n \choose n}$  words of the form (for convenience we write $e_0$ for $t$):
$$e_{i_1}e_{i_1-1}\ldots e_1e_0e_{i_2}e_{i_2-1}\ldots e_1e_0
e_{i_s}e_{i_s-1}\ldots e_1e_0 \cdot$$ 
$$\cdot e_{a_1}e_{a_1-1}\ldots e_{b_1+1}e_{b_1}\ldots
e_{a_k}e_{a_k-1}\ldots e_{b_k+1}e_{b_k}$$
where $s,k \geq 0$, 
$0\leq i_1 \leq i_2\ldots < i_s < a_1 < a_2 <\ldots <  a_k < n$, 
$0 < b_i \leq a_i< n$, and $0< b_1 < b_2 <\ldots <  b_k < n$.

\end{enumerate}
\end{theorem} 

We finish this section by offering the following very useful observation
(compare \cite{P-S-2}).

\begin{proposition}\label{3.14}
Consider a 3-manifold $(M,\{x_i\}_{i=1}^{2n})$ and 
let $x_{2n+1}$ and $x_{2n+2}$ lie on the same boundary component of $M$.
Consider the $R$-homomorphism of RKBSM
$$i_{\#}: {\cal S}_{2,\infty}(M,\{x_i\}_{i=1}^{2n},R,A) \to
{\cal S}_{2,\infty}(M,\{x_i\}_{i=1}^{2n+2},R,A)$$
 generated by the identity map and with convention that $i_{\#}(L)$
has the point $x_{2n+1}$ connected to the point $x_{2n+2}$ 
by a framed arc close to
boundary (we push out of the boundary framed arc joining $x_{2n+1}$ and
$x_{2n+2}$ in $\partial M$).\ \
Then\\
$i_{\#}$ is a monomorphism if one assumes that
$A^2+A^{-2}$ is not an anihilator of any non-zero element of 
${\cal S}_{2,\infty}(M,\{x_i\}_{i=1}^{2n},R,A)$ 
(i.e. $(A^2+A^{-2})x=0 \Rightarrow x=0$).

\end{proposition}
\begin{proof} Consider the $R$-homomorphism 
$$i_{\#}': {\cal S}_{2,\infty}(M,\{x_i\}_{i=1}^{2n+2},R,A) \to
{\cal S}_{2,\infty}(M,\{x_i\}_{i=1}^{2n},R,A)$$ given by connecting
$x_{2n+1}$ and $x_{2n+2}$ in $\partial M$ and pushing it inside $M$.
Now clearly $i_{\#}'i_{\#}(L)=(-A^2-A^{-2})(L)$, thus 
$i_{\#}'i_{\#}(u)=(-A^2-A^{-2})(u)$ for any $u\in 
{\cal S}_{2,\infty}(M,\{x_i\}_{i=1}^{2n},R,A)$. Therefore $i_{\#}'i_{\#}$
is a monomorphism iff $A^2+A^{-2}$ is not an anihilator of any non-zero 
element of ${\cal S}_{2,\infty}(M,\{x_i\}_{i=1}^{2n},R,A)$. 
A monomorphism of $i_{\#}$ follows from a monomorphism of $i_{\#}'i_{\#}$.
\end{proof}
  
\section{Torsion in KBSM}
In all of the  examples above the module is torsion free except in 
the case of $S^1 \times S^2$.
In fact a non-separating $S^2$ in $M$ always yields a torsion in  
${\cal S}_{2,\infty}(M)$. It is enough to use the framing relation to see
a torsion: Let $L$ be a framed link cutting a non-separating $S^2$ exactly
in one point. We can twist $S^2$ twice, twisting also the framing of $L$
twice and then undo this by an isotopy of $M$. Thus $(A^6-1)L=0$ in
${\cal S}_{2,\infty}(M)$. It is less obvious that a separating $S^2$ can
often yield a torsion.
\begin{conjecture}[\cite{Kir}] \label{4.1}
If $M=M_1\# M_2$, where $M_i$ is not equal to $S^3$, possibly with holes,
then ${\cal S}_{2,\infty}(M)$ has a torsion element.
\end{conjecture}

We have proven the  Conjecture 4.1 only partially. In an example in 
\cite{Pr-7} we use the first homology groups of summands of the connected
sum. One can extend it employing
$SL(2,C)$ representations of summands (particularly if summands are
hyperbolic).
\begin{theorem} \label{4.2}
\begin{enumerate}
\item[(a)] (\cite{Pr-7}).
If $M_1$ and $M_2$ have first homology groups that are not 
2-torsion groups, then the conjecture holds.
\item[(b)] If there are representations $\rho_i : \pi_1(M_i) \to SL(2,C)$,
such that the image $\rho_i(\pi_1(M_i))$ is not in the center of $SL(2,C)$
(which is composed of $\pm Id$) then the conjecture holds.
\end{enumerate}
\end{theorem}

In the case of $M$ containing an incompressible torus we are able to
show the following theorem (the first homology group is used in the proof 
of (a) and $SL(2,C)$ representations in the proof of (b)).

\begin{theorem} \label{4.3}
\begin{enumerate}
\item[(a)]
Let $M$ be a 3-manifold allowing embedded non-separating torus.
Then ${\cal S}_{2,\infty}(M)$ has a torsion element.
\item[(b)] Let $M$ be a manifold and $\partial M$ contains a torus, 
$\partial_1 M$. Assume that there is a representations 
$\rho : \pi_1(M) \to SL(2,C)$, such that the image 
$\rho(\pi_1(M))$ is not conjugated to upper triangular subgroup
of matrices (Borel subgroup) in $SL(2,C)$ and $\rho(\pi_1(\partial_1 M))$
is not in the center of $SL(2,C)$, \ \ then 
the double of $M$ along $\partial_1 M$ has a torsion element in its KBSM.
In particular the double of the complement of a hyperbolic knot has
a torsion element in its KBSM; see \cite{Ve}.
%We use the fact that trace (ABA^{-1}B^{-1}) $ is not 2 unless $A$ and $B$
%are simultaneously upper (or lower) triangular and diagonal matrix - not
%in the center can be conjugated out of the diagonal.
\end{enumerate}
\end{theorem}
The relative case similar to Theorem 4.3 is considered in Lemma 4.4(b).

We go back now to our main task of analyzing the relative skein module
of the  product of a surface and the interval.

\begin{lemma}\label{4.4}
\begin{enumerate}
\item[(a)] If $M=S^2 \times I$ and not all $x_i$ are on the same boundary 
component of $M$, then the relative Kauffman bracket skein module of 
${\cal S}_{2,\infty}(M,\{x_i\}_{i=1}^{2n},Z[A^{\pm 1}],A)$ 
has a torsion.
\item[(b)] 
The relative Kauffman bracket skein module of 
$M=T^2 \times I$,\\ 
${\cal S}_{2,\infty}(M,\{x_i\}_{i=1}^{2n},Z[A^{\pm 1}],A)$ 
has a torsion if not all $x_i$ are on the same boundary component of $M$. 
\end{enumerate}
\end{lemma} 
\begin{proof}
\begin{enumerate}
\item[(a-1)]  For $n=1$ and $M=S^2 \times I$ one uses the standard
``Dirac trick" that is if we twist framing twice on the arc joining $x_1$
with $x_2$ in any relative link, we get the relative link ambient isotopic
to $L$ thus $(A^6-1)L=0$. $L$ is not zero as the following easy argument
shows: \ $L$ represents a nontrivial element in the relative homology
group with $Z_2$ coefficients, $H_1(M,\{x_1,x_2\};Z_2)$. 
On the other hand we have
a $Z$-homomomorphism ${\cal S}_{2,\infty}(M,\{x_1,x_2\}) \to
CH_1(M,\{x_1,x_2\};Z_2)$ 
given by sending an unoriented framed link $L$ to an element in
$H_1(M,Z_2)$ it represent and $A$ to $\omega_A$ in $C$, where
$\omega_A$ is the primitive sixth root of unity (that is it satisfies:
$\omega_A + \omega_A^{-1}  = 1$). 
Thus $L \neq 0$ in ${\cal S}_{2,\infty}(M,\{x_1,x_2\})$;. \cite{Pr-1}.
\item[(b-1)] We will describe the construction in details as it has
several uses and generalizations.\\
Let $T_b^2=T^2 \times{1/2}$ be the ``middle" torus and $L$ a relative link
which cuts $T_b^2$ in exactly one point. Further let $\lambda$ be a
noncontractible curve on $T_b^2$. We can represent the link $L\sqcup \lambda$
putting $\lambda$ ``on the top of" $L$ or ``below" $L$. We use the fact that
on the torus $\lambda$ can be isotoped on the ``other" side of $L \cap T_b^2$.
Thus in the RKBSM, $0= (A-A^{-1})({L\lambda}^{-1} - L{\lambda})$ as illustrated
in Fig. 4.1. 
\begin{center}
\centerline{\psfig{figure=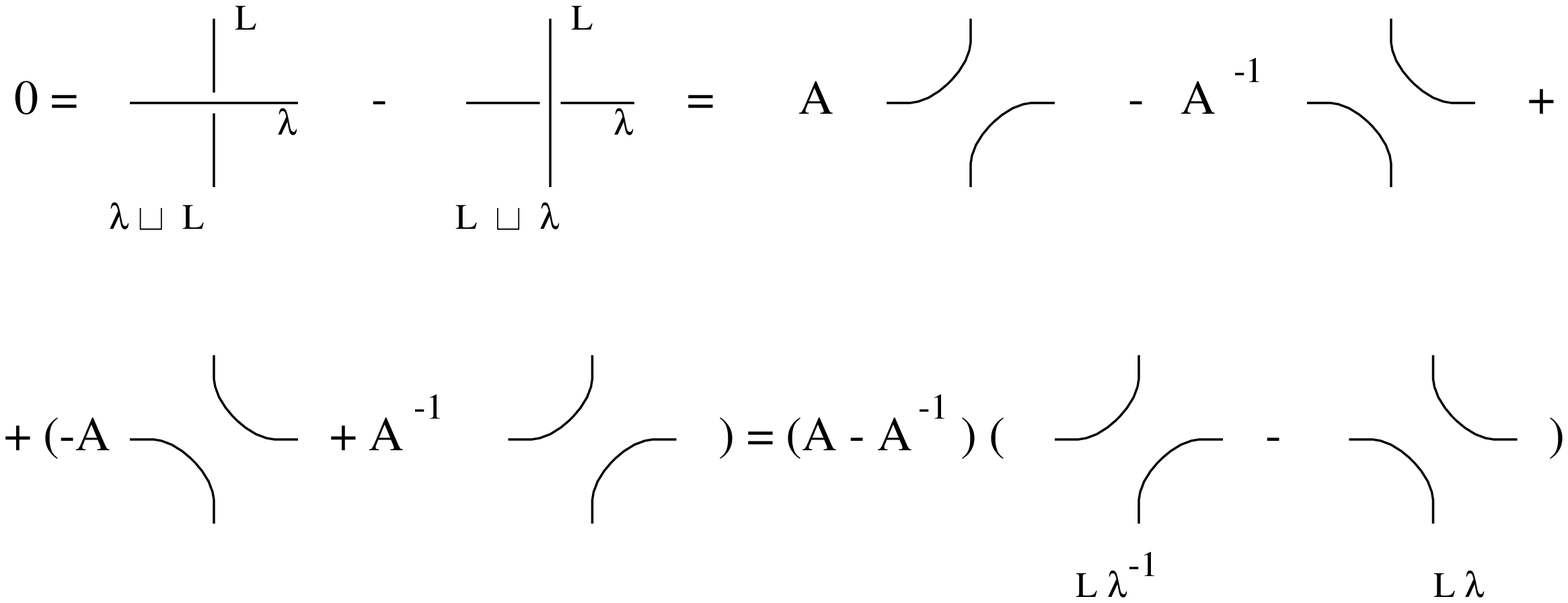,height=4.9cm}}
Figure 4.1
\end{center}

It remain to see that $L\lambda \neq L{\lambda}^{-1}$ in
${\cal S}_{2,\infty}(M,\{x_i\}_{i=1}^{2},Z[A^{\pm 1}],A)$. We show that 
it is not zero even for $A=-1$. 
In fact ${\cal S}_{2,\infty}(M,\{x_i\}_{i=1}^{2} ;Z,-1)$ 
is isomorphic to $Z[Z\oplus Z]$ (it follows from \cite{P-S-1,P-S-2}.
The isomorphism $\phi : {\cal S}_{2,\infty}(M,\{x_i\}_{i=1}^{2} ;Z,-1) \to
Z[Z\oplus Z]$ send an arc, $\gamma$, joining $x_1$ with $x_2$ 
(with an orientation from $x_1$ to $x_2$) to the 
(minus) element in $H_1(M,Z)$ represented by the arc $\gamma$ with endpoints
connected by the straight vertical line joining $x_1$ with $x_2$. 
In fact we only need much weaker fact that
$\phi $ is well defined and $\phi (L\lambda)\neq \phi (L{\lambda}^{-1})$ in
$H_1(M,Z)$. 

\item[(b-2)] Let $x_1, x_2 \in T^2 \times \{1\}$ and 
$x_3, x_4 \in T^2 \times \{0\}$. Consider a link $K$ composed of two arcs
(joining $x_1$ with $x_3$, and $x_2$ with $x_4$, each cutting 
$T_b^2=T^2 \times{1/2}$ in exactly one point ($c_1$ and $c_2$ respectively). 
We choose $\lambda$ like in (b-1) and again use ``two sides" of $c_1$ and $c_2$
(assumed to be closed together) on the torus. We get after some
computation (see Fig. 4.2).
$(A^2-A^{-2})((K_2^{-1}\lambda K_1' \sqcup K_2' K_1^{-1}) -
                 (K_2^{-1}K_1' \sqcup K_2'\lambda K_1^{-1}))=0$.
We have to show that $(K_2^{-1}\lambda K_1' \sqcup K_2' K_1^{-1}) -    
                 (K_2^{-1}K_1' \sqcup K_2'\lambda K_1^{-1})$ is not $0$.
%It is fairly easy exercise, say using $\mod 2$ relative homology.
%But we present here more general, and very useful setting.
%(It can be used to show that the double of a hyperbolic manifold has
%a torsion in the KBSM, compare \cite{Ve}. 
Let us denote our two links by $L_1$ and $L_2$. To distinguish $L_1$
from $L_2$ in
the RKBSM we use Theorem 4.3 (b) (result which depends on  
$SL(2,C)$ representations)
and some topological reasoning (embedding $M$ in bigger manifold and
extending relative links to links without boundary).\\
Let $M_2$ be a manifold obtained from $M=T^2\times [0,1]$ by adding to it
two 1-handles, first with a core joining $x_1$ with $x_2$ and the second
with a core joining $x_3$ with $x_4$. 
Any relative link in $M$ can be extended to a link
in $M_2$ by adding to it the cores of 1-handles. This gives us a
homomorphism:
$$\psi: {\cal S}_{2,\infty}(M,\{x_1,x_2,x_3,x_4\};R,A) \to
       {\cal S}_{2,\infty}(M_2;R,A).$$
We want to show that $\psi (L_1) \neq \psi (L_2)$ in
${\cal S}_{2,\infty}(M_2;Z[A^{\pm 1}],A)$. By Universal Coefficient Property,
it suffice to show this in ${\cal S}_{2,\infty}(M_2;C,-1)$.
Clearly $M_2$ is the double of the manifold $M_1$ obtained from $T^2 \times
[\frac{1}{2},1]$ by adding a 1-handle with a core joining $x_1$ with $x_2$.
The fundamental group $\pi_1(M_1)= (Z\oplus Z)* Z$. We can use Theorem 4.3(b),
because there exists a required representation $\rho : \pi_1(M) \to SL(2,C)$.
We can build the representation concretely, to distinguish 
$\psi (L_1)$ from $\psi (L_2)$. 
For example we can send generators of $(Z\oplus Z)$ to the matrix
\[  \left[ \begin{array}{cc}
1 & 1 \\
0 & 1
   \end{array}
\right]
\]
 
and the generator of the fundamental group composed of $K_1,K_1'$ 
and the core of the added handle to \[  \left[ \begin{array}{cc}
1 & 0 \\
1 & 1
   \end{array}
\right]
\] (compare \cite{Ve}).

\begin{center}
\centerline{\psfig{figure=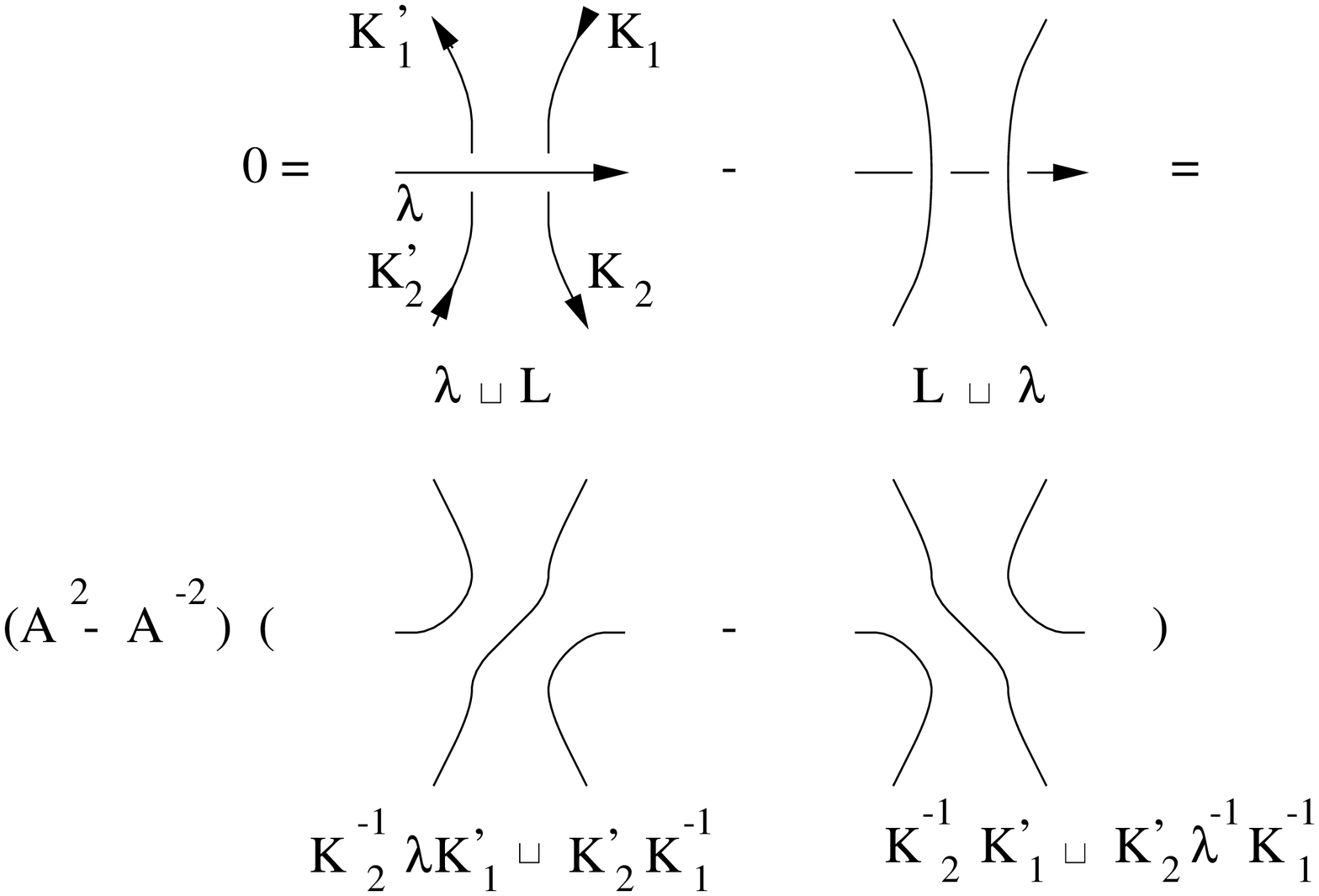,height=7.1cm}}
Figure 4.2
\end{center}

\item[(a-n),(b-n)] For any other number of points, $2n$, Lemma 4.4 
 follows immediately from (a-1) or (b-1) for an odd $n$,
and from (a-2) or (b-2) for an even $n$. 
\end{enumerate} 
\end{proof} 

\begin{problem}\label{4.5}
Let $F$ be an oriented closed surface of genus greater than 1, and not
all points $x_i$ lie on the same boundary component of $F\times [0,1]$.
Does the RKBSM, ${\cal S}_{2,\infty}(F\times [0,1],\{x_i\}_1^{2n})$ has
a torsion element?
\end{problem}
The answer to the above question can shed a light into the role of
incompressible surfaces in the structure of Kauffman bracket
skein modules of 3-manifolds.

\ \\
\ \\
Department of Mathematics\\
George Washington University  \\
Washington, DC 20052 \\
e-mail: przytyck@math.gwu.edu

\end{document}